\documentclass[10pt]{article}
\usepackage{amsmath}
\usepackage{color}
\usepackage{graphicx}
\usepackage{fancyhdr}
\usepackage[english]{babel}
\usepackage[autostyle]{csquotes}
\usepackage{verbatim}
\usepackage{amsfonts,longtable}
\usepackage{epsfig}
\usepackage{mathrsfs}
\usepackage{amsfonts}
\usepackage{color}

\def\R{\hbox{{\rm I}\kern-0.2em{\rm R}\kern0.2em}}

\def\d{{\rm d}}

\def\bn{\begin{equation}}
\def\en{\end{equation}}
\def\bny{\begin{eqnarray}}
\def\eny{\end{eqnarray}}
\def\be{\begin{eqnarray*}}
\def\ee{\end{eqnarray*}}
\def\bc{\begin{center}}
\def\ec{\end{center}}
\def\X{{\cal X}} 
\def\p{\partial}
\def\({\left(}
\def\){\right  )}
\def\[{\left[}
\def\]{\right]}
\def\bc{\begin{center}}
\def\ec{\end{center}}

\newtheorem{dfn}{Definition}[section]
\newtheorem{thm}{Theorem}[section]
\newtheorem{rem}{Remark}[section]
\newtheorem{pro}{Proposition}[section]

\newtheorem{cor}{Corollary}[section]
\newtheorem{lem}{Lemma}[section]
\newtheorem{exm}{Example}[section]

\def\bn{\begin{equation}}
\def\en{\end{equation}}
\def\bny{\begin{eqnarray}}
\def\eny{\end{eqnarray}}
\def\ba{\begin{array}{lllll}}
\def\ea{\end{array}}
\def\be{\begin{eqnarray*}}
\def\ee{\end{eqnarray*}}
\def\bdn{\begin{dfn}}
\def\edn{\end{dfn}}
\def\btm{\begin{thm}}
\def\etm{\end{thm}}
\def\bpf{\begin{proof}}
\def\epf{\end{proof}}
\def\bpn{\begin{pro}}
\def\epn{\end{pro}}
\def\brk{\begin{rem}}
\def\erk{\end{rem}}
\def\bcy{\begin{cor}}
\def\ecy{\end{cor}}
\def\blm{\begin{lem}}\def\elm{\end{lem}}
\def\bex{\begin{exm}}
\def\eex{\end{exm}}

\def\X{{\cal X}}  
  
 \def\R{\mathscr{R}}
\def\S{{\cal S}}

\begin{document}
\bc \Large{{\bf  Invariance analysis, reductions and conservation laws for second-order systems of ordinary difference equations.

}}\ec
\medskip

\bc

J J H Bashingwa$^{1,a}$ and A H Kara$^{1,2,b}$\\
$^{1}$School of Mathematics, University of the Witwatersrand, \\
Johannesburg, Wits 2001, South Africa.\\

%
%
$^{2}$Department of Mathematics and Statistics, King Fahd University of
Petroleum and Minerals, Dhahran 31261, Saudi Arabia\\

{\small $^a$ jean.bashingwa@wits.ac.za, $^b$ Abdul.Kara@wits.ac.za}
\ec
 
\begin{abstract}
We present a method to obtain symmetries for second-order systems of ordinary difference equations and how to use them to reduce the order. We also introduce  a technique of finding conservation laws for such systems.
\end{abstract}
\section{Introduction}
Physical laws are described by differential equations (DEs). Lie group theory provides us with powerful tools for obtaining analytical solutions of such equations \cite{Olver93}. Over the last 30 years, a considerable amount of work has been invested into applying Lie's theory to solve and classify difference  equations ($\Delta$Es) (see \cite{maeda,levi,qui,hyd,vlad} and references therein). 

The use of symmetry methods for ordinary difference equations  (O$\Delta$Es ) has been introduced by Maeda \cite{maeda}.  He showed that the resulting linearized symmetry condition (LSC) amounts to a set of functional equations which is hard to solve in general. In \cite{hyd}, Hydon introduced a technique to solve LSC and  obtain symmetries in closed form by  repetitive differentiations. Now we can say that the symmetry analysis for $\Delta$Es is well documented. For  second-order  O$\Delta$Es a full classification according to their point symmetries exist in the literature \cite{vlad}. 

It has been proved in \cite{hyd} that every second-order linear homogeneous  O$\Delta$E has an eight-dimension Lie algebra isomorphic to $\mathfrak{s}\mathfrak{l}(3)$. However, for second-order system of difference equations (S$\Delta$Es) things may be differents. We shall prove this in Section 3


\section{Groundwork}

Let us consider an  $N$-th   order system of $r$ $\Delta$Es

\bn
 x_{n+N}^i=\omega_i(n,x_n ^1,...,x_n^r,x_{n+1}^1,...,
 x_{n+1}^r,...,
 x_{n+N-1}^1,...,x_{n+N-1}^r),\quad i=1,...,r.\label{eq3}\en
 
 We assume that for each $\omega_i$ there exist at least one $x_n^j$
($i,j=1,...,r$)   such that $\frac{\p \omega_i	}{\p x_n^j}\neq 0 $.

Consider a point transformation 

$$ \Gamma_\epsilon :X \mapsto \hat{X}(X,\epsilon ) 
$$

where $X =( x_n^ 1 ,  . . . ,x_n^r)$ are continuous variables. $\Gamma$ will be called 
 one-parameter Lie group of transformations if it satisfies the following properties:
 
\begin{itemize}
\item[$\bullet$] $\Gamma_0$ is the identity map, i.e, $\hat{X}=X$ for $\epsilon=0.$
\item[$\bullet$] $\Gamma_\mu \Gamma_\nu =\Gamma_{\mu+ \nu}$ for every $\mu$ and $\nu$ close to 0.

\item[$\bullet$] Each $\hat{x}_n^i$ can be expanded as a Taylor series in a neighbourhood of $\epsilon=0$

\end{itemize}

Therefore, we have

\bn  \begin{array}{lll}
 \hat{x}_{n+j}^i=x_{n+j}^i + \epsilon  \S^jQ_i(n,x_n ^1,...,x_n^r,x_{n+1}^1,...,
 x_{n+1}^r,...,
 x_{n+N-1}^1,...,x_{n+N-1}^r)+O(\epsilon^2)\end{array} \en

where  $Q_i$ are continuous functions which we shall refer to as  characteristics, $i=1,...,r$ ,$j=1,...,N$ and $\S$ is the \lq\lq shift\rq\rq operator.
 It is defined as follow

\bn  \S: n \mapsto n+1 ,\quad  
\S^k (x_n^i) =x_{n+k}^i   \en

The symmetry condition for the S$\Delta$Es \eqref{eq3} is

\bn  
 \hat{x}_{n+N}^i=\omega_i({n},\hat{x}_n ^1,...,\hat{x}_n^r,\hat{x}_{n+1}^1,...,
 \hat{x}_{n+1}^r,...,
 \hat{x}_{n+N-1}^1,...,\hat{x}_{n+N-1}^r), \quad i=1,...,r.\label{eq6}\en
whenever \eqref{eq3} holds.

Lie symmetries are obtained by linearizing the symmetry condition \eqref{eq6} about the identity. We have the following system of linearized symmetry condition (SLSC)

\bn 
  \S^N (Q_i) -X\omega_i=0, \quad i=1,...,r.\label{eq4}\en 
 where the symmetry generator $X$ is given by
 
 \bn X= \sum\limits_{j=0}^{N-1}\left(\sum\limits_{i=1}^{r}\S^j(Q_i)
 \frac{\p}{\p x_{n+j}^i}\right)\en
 
\bdn A function $w_n$  is invariant function under the Lie group of transformations $\Gamma$ if \edn

\bn X(w_n)=0,\en

where $w_n$ can be found by solving the characteristic equation

\bn\frac{\d x_n^1}{Q_1} =...=\frac{\d x_{n}^r}{Q_r}=\frac{\d x_{n+1}^1}{\S (Q_1)}=...=\frac{\d x_{n+1}^r}{\S (Q_r)}=....=\frac{\d x_{n+N-1}^1}{\S^{N-1}(Q_1)}=...=\frac{\d x_{n+N-1}^r}{\S^{N-1}(Q_r)}=\frac{w_n}{0} \en

The first integral for the system  \eqref{eq3} is given by

\bn (\S - Id)\phi(n,x_n,y_n,x_{n+1},y_{n+1})=0 \label{eq7} \en

whenever \eqref{eq3} holds.

I Section 3, we shall use the condition \eqref{eq7} to develop a constructive technique for obtaining first integrals.

\brk \erk In this paper we shall consider Lie  point symmetry, i.e,    the characteristics are given by      $Q_i(n,x_n^1,...,x_n^r)$

We  refer the reader  to \cite{Olver93} for more informations on symmetry methods
for differential equations.
\section{ symmetries and reductions} 
Consider a second-order system of 2 $\Delta$Es

\bn\ba  x_{n+2}=\omega_1(n,x_n,y_n,x_{n+1},y_{n+1}),\\
 y_{n+2}=\omega_2(n,x_n,y_n,x_{n+1},y_{n+1}).\ea \en

We assume that $\frac{\p \omega_1}{\p x_n} \neq 0$ or $\frac{\p \omega_1}{\p y_n} \neq 0$ and $\frac{\p \omega_2}{\p x_n} \neq 0$ or $\frac{\p \omega_2}{\p y_n} \neq 0$, so the system is of second order.

The SLCS \eqref{eq4} reduces to 
\bny   \S^2( Q_1)-Q_1\omega_{1,x_n}-Q_2\omega_{1,y_n} -\S (Q_1)\omega_{1,x_{n+1}}-\S (Q_2)\omega_{1,y_{n+1}}=0,\label{eq8}\\ 
\S^2 (Q_2)-Q_1\omega_{2,x_n}-Q_2\omega_{2,y_n} -\S (Q_1)\omega_{2,x_{n+1}}-\S (Q_2)\omega_{2,y_{n+1}}=0. \label{eq9}\eny

where $g_{,x}=\frac{\p g}{\p x}, Q_1(n,x_n,y_n)$ and $Q_2(n,x_n,y_n)$

The functional equations \eqref{eq8} \& \eqref{eq9} contain  functions $Q_1$ and $Q_2$ with different pairs of arguments making them difficult to solve. 
For concreteness, if for instance the discrete variable $n$ stands for \lq \lq state\rq\rq,  $Q_1(n,x_n,y_n)$ and $\S(Q_1)\equiv Q_1(n+1,x_{n+1},y_{n+1})$ belong to two different states. 

To overcome this, we proceed as follows

$\bullet$ Step1:  elimination of  $\S^2(Q_1)$ and $\S^2(Q_2)$

 We differentiate (total differentiation) \eqref{eq8} and  \eqref{eq9} with respect to $x_n$ and $y_n$ respectively, keeping $\omega_1$ and $\omega_2$ fixed . Here, we take  $x_{n+1}$  as function of $x_n,y_n,y_{n+1},\omega_1,\omega_2$ and 
$y_{n+1}$  as function of $x_n,y_n,x_{n+1},\omega_1,\omega_2$
 
 The total derivative operators are given by
 \bn\ba \frac{\d}{\d x_n}=\frac{\p}{\p x_n}+\frac{\p x_{n+1}}{\p x_n}	\frac{\p}{\p x_{n+1}} + \frac{\p y_{n+1}}{\p x_n}	\frac{\p}{\p y_{n+1}} +...\\
 \frac{\d}{\d y_n}=\frac{\p}{\p y_n}+\frac{\p x_{n+1}}{\p y_n}	\frac{\p}{\p x_{n+1}} + \frac{\p y_{n+1}}{\p y_n}	\frac{\p}{\p y_{n+1}} +...
 \ea\en
 
 In this case, this simplifies to 
 \bny  \frac{\d}{\d x_n} =\frac{\p}{\p x_n} - \left(\frac{\omega_{1,x_n}}{\omega_{1,x_{n+1}}} +\frac{\omega_{2,x_n}}{\omega_{2,x_{n+1}}}  \right)\frac{\p}{\p x_{n+1}} -
 \left(\frac{\omega_{1,x_n}}{\omega_{1,y_{n+1}}} +\frac{\omega_{2,x_n}}{\omega_{2,y_{n+1}}}  \right)\frac{\p}{\p y_{n+1}}\label{eq10}\\ 
 \frac{\d}{\d y_n} =\frac{\p}{\p y_n} - \left(\frac{\omega_{1,y_n}}{\omega_{1,x_{n+1}}} +\frac{\omega_{2,y_n}}{\omega_{2,x_{n+1}}}  \right)\frac{\p}{\p x_{n+1}} -
 \left(\frac{\omega_{1,y_n}}{\omega_{1,y_{n+1}}} +\frac{\omega_{2,y_n}}{\omega_{2,y_{n+1}}}  \right)\frac{\p}{\p y_{n+1}}
 \label{eq11}\eny

So we apply the operator \eqref{eq10} to \eqref{eq8} and \eqref{eq11} to \eqref{eq9} keeping $\omega_1$ and $\omega_2$ fixed . This leads to the determining system

\bn\ba  [ Q_1\omega_{1,x_n}+Q_2\omega_{1,y_n}]_{,x_n} 

+\S (Q_1)\omega_{1,x_{n+1},x_n}+\S (Q_2)\omega_{1,y_{n+1},x_n}

-\left(\frac{\omega_{1,x_n}}{\omega_{1,x_{n+1}}} +\frac{\omega_{2,x_n}}{\omega_{2,x_{n+1}}}  \right)[ Q_1\omega_{1,x_n}+Q_2\omega_{1,y_n}]_{,x_{n+1}}\\

-\left(\frac{\omega_{1,x_n}}{\omega_{1,x_{n+1}}} +\frac{\omega_{2,x_n}}{\omega_{2,x_{n+1}}}  \right)
[\S (Q_1)\omega_{1,x_{n+1},x_{n+1}}+\S (Q_2)\omega_{1,y_{n+1},x_{n+1}}  
+[\S (Q_1)]_{,x_{n+1}}\omega_{1,x_{n+1}}+[\S (Q_2)]_{,x_{n+1}}\omega_{1,y_{n+1}}]\\

-\left(\frac{\omega_{1,x_n}}{\omega_{1,y_{n+1}}} +\frac{\omega_{2,x_n}}{\omega_{2,y_{n+1}}}  \right)
[ Q_1\omega_{1,x_n}+Q_2\omega_{1,y_n}]_{y_{n+1}}

-\left(\frac{\omega_{1,x_n}}{\omega_{1,y_{n+1}}} +\frac{\omega_{2,x_n}}{\omega_{2,y_{n+1}}}  \right)[\S (Q_1)\omega_{1,x_{n+1},y_{n+1}}+\S (Q_2)\omega_{1,y_{n+1},y_{n+1}}  
\\+[\S (Q_1)]_{,y_{n+1}}\omega_{1,x_{n+1}}+[\S (Q_2)]_{,y_{n+1}}\omega_{1,y_{n+1}}]
=0,\label{eq12}\ea\en

\bn\ba
[ Q_1\omega_{2,x_n}+Q_2\omega_{2,y_n}]_{,y_n} 

+\S (Q_1)\omega_{2,x_{n+1},y_n}+\S (Q_2)\omega_{2,y_{n+1},y_n}

-\left(\frac{\omega_{1,y_n}}{\omega_{1,x_{n+1}}} +\frac{\omega_{2,y_n}}{\omega_{2,x_{n+1}}}  \right)[ Q_1\omega_{2,x_n}+Q_2\omega_{2,y_n}]_{,x_{n+1}}\\

-\left(\frac{\omega_{1,y_n}}{\omega_{1,x_{n+1}}} +\frac{\omega_{2,y_n}}{\omega_{2,x_{n+1}}}  \right)
[\S (Q_1)\omega_{2,x_{n+1},x_{n+1}}+\S (Q_2)\omega_{2,y_{n+1},x_{n+1}}  
+[\S (Q_1)]_{,x_{n+1}}\omega_{2,x_{n+1}}+[\S (Q_2)]_{,x_{n+1}}\omega_{2,y_{n+1}}]\\

-\left(\frac{\omega_{1,y_n}}{\omega_{1,y_{n+1}}} +\frac{\omega_{2,y_n}}{\omega_{2,y_{n+1}}}  \right)
[ Q_1\omega_{2,x_n}+Q_2\omega_{2,y_n}]_{y_{n+1}}

-\left(\frac{\omega_{1,y_n}}{\omega_{1,y_{n+1}}} +\frac{\omega_{2,y_n}}{\omega_{2,y_{n+1}}}  \right)[\S (Q_1)\omega_{2,x_{n+1},y_{n+1}}+\S (Q_2)\omega_{2,y_{n+1},y_{n+1}}  
\\+[\S (Q_1)]_{,y_{n+1}}\omega_{2,x_{n+1}}+[\S (Q_2)]_{,y_{n+1}}\omega_{2,y_{n+1}}]
=0.\label{eq13}
 \ea\en

 $\bullet$ Step 2: elimination of  $\S(Q_1)$ and $\S(Q_2)$

We now differentiate \eqref{eq12} \& \eqref{eq13} with respect to $x_n$ and $y_n$ respectively keeping $x_{n+1}$ and $y_{n+1}$ fixed. This means that we apply the operator $\frac{\p}{\p x_n}$ on \eqref{eq12} and $\frac{\p}{\p y_n}$ on \eqref{eq13}.  For a second-order S$\Delta$Es we need at most to differentiate four times. After separating with respect to $x_{n+1}$ and $y_{n+1}$  the resulting equations, we obtain a system of 
DEs in $Q_1$ and $Q_2$ which is solvable by hand or by using a computer algebra package.

$\bullet$ Step 3: explicit form of constant of integration

When integrating in step2 to obtain the characteristics $Q_1$ and $Q_2$, we have constant of integration which appear to be functions of $n$. To obtain their explicit form, we need to substitute the results obtained in step2 in \eqref{eq12} \& \eqref{eq13}.

\section{Applications}
\subsection{Example1}
Consider the most general homogeneous  2nd  order linear  system of difference equations 

\bn \ba x_{n+2}= a_1(n)x_n + a_2(n)y_n +a_3(n)x_{n+1}+a_4(n)y_{n+1},\\
y_{n+2}= b_1(n)x_n + b_2(n)y_n +b_3(n)x_{n+1}+b_4(n)y_{n+1}\label{eqc1} \ea\en where $a_i(n),b_i(n), \quad i=1,...,4$ are arbitrary functions.

 One can readily verify that  the determining system \eqref{eq12} \& \eqref{eq13} amounts to
 \bn Q_{1,x_nx_n}=Q_{2,x_n x_n}=0, \quad Q_{1,y_ny_n}=Q_{2,y_ny_n}=0\en. 
 
 Therefore
 
 \bn Q_1(n,x_n,y_n)=C_1 x_n +C_2y_n +F_1(n), \quad Q_2(n,x_n,y_n)=C_3 x_n +C_4y_n +F_2(n) \label{eqc3}\en 
 where $C_i,\quad i=1...4$ are constants.

 The characteristics in \eqref{eqc3} must satisfy the SLSC \eqref{eq8} \& \eqref{eq9}. Hence, we have

 \bn\ba
 F_1(n+2)-[a_1(n)F_1(n)+a_2(n)F_2(n)+a_3(n)F_1(n+1)+a_4(n)F_2(n+1)]=0,\\
 F_2(n+2)-[b_1(n)F_1(n)+b_2(n)F_2(n)+b_3(n)F_1(n+1)+b_4(n)F_2(n+1)]=0. \label{eqc4}\ea \en
 
 and $$  C_1=C_4, \qquad C_2=C_3=0 $$
 So, \eqref{eqc3} simplifies to 
 
  \bn Q_1(n,x_n,y_n)=C_1 x_n  +F_1(n), \quad Q_2(n,x_n,y_n)=C_1y_n +F_2(n) \label{eqc5}\en 
 The first generator of symmetry for a second-order homogeneous linear system \eqref{eqc4} is the scaling symmetry  given by
 
 $$X=x_n\p _{x_n}+ y_n \p_{ y_n}$$ 
 The system \eqref{eqc4}, which governs the remaining generators  of the Lie point symmetry for the system \eqref{eqc1}, is of  second order in $F_1$ and $F_2$. Its general solution is 
 
 \bn F_1(n)=g_1(n,K_1,K_2,K_3,K_4), \quad F_2(n)=g_2(n,K_1,K_2,K_3,K_4) \en where $K_1,\quad i=1,...,4$ are constants.
 
 So, the most large Lie algebra of symmetry generators which can be obtained from a homogeneous second-order system of 2 difference equations has dimension five.
 
For clarification let us consider 

 $a_1(n)=a_3(n)=a_4(n)=0, \quad a_2(n)=1$ and $ b_2(n)=b_3(n)=b_4(n)=0,\quad b_1(n)=1$. The system \eqref{eqc1} becomes \bn x_{n+2}=y_n, \quad y_{n+2}=x_n \en 
 
 The system which governs the remaining generators of the   Lie point symmetry in this case is given by
 
 \bn F_1(n+2)-F_2(n)=0, \quad F_2(n+2)-F_1(n)=0\en The general solutions for this system will be 
 
 \bn\ba F_1(n)=&\frac{C_1[1+(-1)^n+i^n+(-i)^n]}{4}+\frac{C_2[1-(-1)^n-i(i)^n+i(-i)^n]}{4}+\frac{C_3[1+(-1)^n-i^n-(-i)^n]}{4} \\
 &\frac{C_4[1-(-1)^n+i(i)^n-i(-i)^n]}{4}\\
 
 F_2(n)=&\frac{C_3[1+(-1)^n+i^n+(-i)^n]}{4}+\frac{C_4[1-(-1)^n-i(i)^n+i(-i)^n]}{4}+\frac{C_1[1+(-1)^n-i^n-(-i)^n]}{4} \\
 &\frac{C_2[1-(-1)^n+i(i)^n-i(-i)^n]}{4}\ea\en
 Therefore we have 5 generators of the Lie point symmetry spanned by
 
 \bn\ba \X_0 =x_n{\p _{x_n}}+ y_n \p_{ y_n}\\ \X_1=\frac{[1+(-1)^n+i^n+(-i)^n]}{4} \p_{x_n}+ \frac{[1+(-1)^n-i^n-(-i)^n]}{4} \p_{y_n}\\
 
 \X_2=\frac{[1-(-1)^n-i(i)^n+i(-i)^n]}{4}\p_{x_n}+\frac{[1-(-1)^n+i(i)^n-i(-i)^n]}{4}\p_{y_n}\\
 \X_3=\frac{[1+(-1)^n-i^n-(-i)^n]}{4}   \p_{x_n} +\frac{[1+(-1)^n+i^n+(-i)^n]}{4}\p_{y_n}\\
 \X_4= \frac{[1-(-1)^n+i(i)^n-i(-i)^n]}{4}  \p_{x_n} +\frac{[1-(-1)^n-i(i)^n+i(-i)^n]}{4}\p_{y_n}\ea \en

\subsection{Example 2}

Consider the system
\bn\ba x_{n+2}=\frac{x_{n}y_{n+1}+1}{x_n +y_{n+1}}\\

y_{n+2}=\frac{y_{n}x_{n+1}+1}{y_n +x_{n+1}} \label{eqb0}\ea \en
From the invariants 

The determining system give 

\bny 
-Q_{{2,y_{{n}}}}{x_{{n+1}}}^{2}{y_{{n}}}^{2}+{\it SQ}_{{1,x_{{n+1}}}}{
x_{{n+1}}}^{2}{y_{{n}}}^{2}+2\,Q_{{2}}{x_{{n+1}}}^{2}y_{{n}}-2\,{\it 
SQ}_{{1}}{y_{{n}}}^{2}x_{{n+1}}+Q_{{2,y_{{n}}}}{x_{{n+1}}}^{2}+Q_{{2,y
_{{n}}}}{y_{{n}}}^{2}-{\it SQ}_{{1,x_{{n+1}}}}{x_{{n+1}}}^{2}\nonumber\\
-{\it SQ}
_{{1,x_{{n+1}}}}{y_{{n}}}^{2}-2\,Q_{{2}}y_{{n}}+2\,{\it SQ}_{{1}}x_{{n
+1}}-Q_{{2,y_{{n}}}}+{\it SQ}_{{1,x_{{n+1}}}}=0\label{eqb1}\\
-Q_{{1,x_{{n}}}}{x_{{n}}}^{2}{y_{{n+1}}}^{2}+{\it SQ}_{{2,y_{{n+1}}}}{
x_{{n}}}^{2}{y_{{n+1}}}^{2}+2\,Q_{{1}}{y_{{n+1}}}^{2}x_{{n}}-2\,{\it 
SQ}_{{2}}{x_{{n}}}^{2}y_{{n+1}}+Q_{{1,x_{{n}}}}{x_{{n}}}^{2}+Q_{{1,x_{
{n}}}}{y_{{n+1}}}^{2}-{\it SQ}_{{2,y_{{n+1}}}}{x_{{n}}}^{2}\nonumber\\
-{\it SQ}_{
{2,y_{{n+1}}}}{y_{{n+1}}}^{2}-2\,Q_{{1}}x_{{n}}+2\,{\it SQ}_{{2}}y_{{n
+1}}-Q_{{1,x_{{n}}}}+{\it SQ}_{{2,y_{{n+1}}}}=0\label{eqb2}
\eny

Differentiating twice \eqref{eqb1} with respect to $x_n$ and twice \eqref{eqb2} with respect to $y_n$ keeping $x_{n+1}$ and $y_{n+1}$ fixed we obtain, after separating with respect to  $x_{n+1}$ and $y_{n+1}$ the following system of DEs

\bn \ba  Q''_{{1}}+x_{{n}}{}Q'''_{{1}}+2\,{\frac {Q_{{1}}}{{x_{{n}
}}^{2}}}-2\,{\frac {Q'_{{1}}}{x_{{n}}}}+
{\frac {Q''_{{1}}}{{x_{{n}}}^
{2}}}-{\frac {Q'''_{{1}}}{x_{{
n}}}}=0\\

Q''_{{2}}+y_{{n}}Q'''_{{2}}+2\,{\frac {Q_{{2}}}{{y_{{n}
}}^{2}}}-2\,{\frac {Q'_{{2}}}{y_{{n}}}}+
{\frac {Q''_{{2}}}{{y_{{n}}}^
{2}}}-{\frac {Q'''_{{2}}}{y_{{
n}}}}=0
 \ea\en

whose most general solutions are  

\bn\ba  Q_1(n,x_n)= F_1(n)x_n +F_2(n)(x_n^2-1)\ln \frac{x_n+1}{x_n-1}+F_3(n)(x_n^2-1)\\
Q_2(n,y_n)= F_4(n)y_n +F_5(n)(y_n^2-1)\ln \frac{y_n+1}{y_n-1}+F_6(n)(y_n^2-1)\label{eqb3}\ea\en

To obtain the nature of functions $F_1,...,F_6$ we substitute  \eqref{eqb3} in \eqref{eqb1} and \eqref{eqb2}. After separating with respect to $x_{n},x_{n+1},y_n$ and $y_{n+1}$ we get the following S$\Delta$Es

\bn \ba-4\,{ F2} ( n ) -2\,{ F4} ( n+1 ) +4\,{
 F5} ( n+1 ) -2\,{ F1} ( n )=0\\

4\,{ F2} ( n ) -2\,{ F4} ( n+1 ) -4\,{ 
F5} ( n+1 ) +2\,{ F1} ( n ) =0\\

-4\,{ F5} ( n ) -2\,{ F1} ( n+1 ) +4\,{
 F2} ( n+1 ) -2\,{ F4} ( n ) =0
\\
4\,{ F5} ( n ) -2\,{ F1} ( n+1 ) -4\,{
F2} ( n+1 ) +2\,{ F4} ( n ) =0
 \ea\en

whose solutions are

\bn \ba F_1(n)=F_4(n)=0\\
F_2(n)=C_1+(-1)^nC_2 \\F_5(n)=C_1-(-1)^nC_2 \label{eqb4} \ea\en

The remaining unknown functions $F_3(n)$ and $F_6(n)$ are determined by substituting \eqref{eqb4} \& \eqref{eqb3} into the SLSC ..... This leads to the S$\Delta$Es 

\bn \ba F_3(n)-F_3(n+2)+F_6(n+1)=0\\F_6(n)+F_3(n+1)-F_6(n+2)=0 \label{eqb5}\ea\en

The general solutions to \eqref{eqb5} is given by

\bn\ba F3(n)=& \frac{1+(-1)^n}{\sqrt{5}}\left\lbrace  
\frac{[(-1+\sqrt{5})^{n-1}+(1+\sqrt{5})^{n-1}]C_3}{2^n}

-\frac{[(-1+\sqrt{5})^{n}-(1+\sqrt{5})^{n}]C_6}{2^{n+1}}\right\rbrace\\

& \frac{-1+(-1)^n}{\sqrt{5}}\left\lbrace  
\frac{[(-1+\sqrt{5})^{n-1}-(1+\sqrt{5})^{n-1}]C_5}{2^n}

-\frac{[(-1+\sqrt{5})^{n}+(1+\sqrt{5})^{n}]C_4}{2^{n+1}} \right\rbrace \\

F6(n)=& \frac{1+(-1)^n}{\sqrt{5}}\left\lbrace  
\frac{[(-1+\sqrt{5})^{n-1}+(1+\sqrt{5})^{n-1}]C_5}{2^n}

-\frac{[(-1+\sqrt{5})^{n}-(1+\sqrt{5})^{n}]C_4}{2^{n+1}}\right\rbrace\\

& \frac{-1+(-1)^n}{\sqrt{5}}\left\lbrace  
\frac{[(-1+\sqrt{5})^{n-1}-(1+\sqrt{5})^{n-1}]C_3}{2^n}

-\frac{[(-1+\sqrt{5})^{n}+(1+\sqrt{5})^{n}]C_6}{2^{n+1}} \right\rbrace  \ea \en

where $C_1,...,C_6$ are arbitrary constants. It follows that the characteristics are given by

\bn \ba Q_1=&(C_1+(-1)^nC_2)(x_n^2-1)\ln \frac{x_n+1}{x_n-1}+ 

\left\lbrace \frac{1+(-1)^n}{\sqrt{5}}\left[ 
\frac{[(-1+\sqrt{5})^{n-1}+(1+\sqrt{5})^{n-1}]C_3}{2^n}

-\frac{[(-1+\sqrt{5})^{n}-(1+\sqrt{5})^{n}]C_6}{2^{n+1}}\right]\right.\\ 

&\left. \frac{-1+(-1)^n}{\sqrt{5}}\left[ 
\frac{[(-1+\sqrt{5})^{n-1}-(1+\sqrt{5})^{n-1}]C_5}{2^n}

-\frac{[(-1+\sqrt{5})^{n}+(1+\sqrt{5})^{n}]C_4}{2^{n+1}} \right]\right\rbrace (x_n^2-1) \\

Q_2=&(C_1-(-1)^nC_2)(y_n^2-1)\ln \frac{y_n+1}{y_n-1}+

\left\lbrace \frac{1+(-1)^n}{\sqrt{5}}\left[  
\frac{[(-1+\sqrt{5})^{n-1}+(1+\sqrt{5})^{n-1}]C_5}{2^n}
-\frac{[(-1+\sqrt{5})^{n}-(1+\sqrt{5})^{n}]C_4}{2^{n+1}}\right]\right.\\

& \left. \frac{-1+(-1)^n}{\sqrt{5}}\left[ 
\frac{[(-1+\sqrt{5})^{n-1}-(1+\sqrt{5})^{n-1}]C_3}{2^n}
-\frac{[(-1+\sqrt{5})^{n}+(1+\sqrt{5})^{n}]C_6}{2^{n+1}} \right] \right\rbrace (y_n^2-1)

\ea\en
Therefore, we have six generators of Lie point symmetry

\bn \ba \X_1=(x_n^2-1)\ln \frac{x_n+1}{x_n-1} \p_{x_n} + (y_n^2-1)\ln \frac{y_n+1}{y_n-1}\p_{y_n} \\

\X_2= (-1)^n(x_n^2-1)\ln \frac{x_n+1}{x_n-1}\p_{x_n}  -(-1)^n(y_n^2-1)\ln \frac{y_n+1}{y_n-1}\p_{y_n}\\

\X_3=   \frac{1+(-1)^n}{\sqrt{5}}\left[ 
\frac{[(-1+\sqrt{5})^{n-1}+(1+\sqrt{5})^{n-1}]}{2^n}

\right]
 (x_n^2-1)\p_{x_n} +\frac{-1+(-1)^n}{\sqrt{5}}\left[ 
\frac{[(-1+\sqrt{5})^{n-1}-(1+\sqrt{5})^{n-1}]}{2^n}
 \right]  (y_n^2-1)
 \p_{y_n}\\

\X_4=    \frac{-1+(-1)^n}{\sqrt{5}}\left[ 

\frac{[(-1+\sqrt{5})^{n}+(1+\sqrt{5})^{n}]}{2^{n+1}} \right] (x_n^2-1)\p_{x_n} + \frac{1+(-1)^n}{\sqrt{5}}\left[  
\frac{[(-1+\sqrt{5})^{n}-(1+\sqrt{5})^{n}]}{2^{n+1}}\right] (y_n^2-1)
 \p_{y_n}\\

\X_5=    \frac{-1+(-1)^n}{\sqrt{5}}\left[ 
\frac{[(-1+\sqrt{5})^{n-1}-(1+\sqrt{5})^{n-1}]}{2^n}

 \right](x_n^2-1)\p_{x_n} + \frac{1+(-1)^n}{\sqrt{5}}\left[  
\frac{[(-1+\sqrt{5})^{n-1}+(1+\sqrt{5})^{n-1}]}{2^n}
\right] (y_n^2-1)
 \p_{y_n}\\

\X_6=    \frac{1+(-1)^n}{\sqrt{5}}\left[ 

\frac{[(-1+\sqrt{5})^{n}-(1+\sqrt{5})^{n}]}{2^{n+1}}\right]  (x_n^2-1)\p_{x_n} +\frac{-1+(-1)^n}{\sqrt{5}}\left[ 
\frac{[(-1+\sqrt{5})^{n}+(1+\sqrt{5})^{n}]}{2^{n+1}} \right] (y_n^2-1)
 \p_{y_n}\label{eq6b}

\ea\en

Each generator in \eqref{eq6b} can be used to reduce the order of \eqref{eqb0}.

Let us consider $\X_1$. By the characteristic method for Partial Differential Equations, the invariants are given by following equation

$$\frac{\d x_n}{(x_n^2-1)\ln \frac{x_n+1}{x_n-1}}=\frac{\d y_n}{(y_n^2-1)\ln \frac{y_n+1}{y_n-1}}=\frac{\d x_{n+1}}{(x_{n+1}^2-1)\ln \frac{x_{n+1}+1}{x_{n+1}-1}}=\frac{\d y_{n+1}}{(y_{n+1}^2-1)\ln \frac{y_{n+1}+1}{y_{n+1}-1}}=\frac{V_n}{0}$$

We get 

$$C_1=\frac{\ln \frac{x_n +1}{x_n -1}}{\ln \frac{y_{n+1}+1}{y_{n+1}-1}},\quad  C_2=\frac{\ln \frac{x_n +1}{x_n -1}}{\ln \frac{y_{n}+1}{y_{n}-1}},\quad  C_3=\frac{\ln \frac{x_n +1}{x_n -1}}{\ln \frac{x_{n+1}+1}{x_{n+1}-1}},\quad V_n=f(C_1,C_2,C_3)$$

where $C_1,C_2,C_3$ are constants.

If we choose $f(C_1,C_2,C_3)=C_1$, we have

\bn u_n=\frac{\ln \frac{x_n +1}{x_n -1}}{\ln \frac{y_{n+1}+1}{y_{n+1}-1}}\label{eq7b} \en

and if  we choose $f(C_1,C_2,C_3)=\frac{C_3}{C_2}$, we have

\bn v_n=\frac{\ln \frac{y_n +1}{y_n -1}}{\ln \frac{x_{n+1}+1}{x_{n+1}-1}}\label{eq8b} \en
%
From \eqref{eq6b} \& \eqref{eq7b} 

we deduce 
\bn\ba  u_{n+1}=
\frac{\ln \frac{x_{n+1}+1}{x_{n+1}-1}}{\ln \frac{(x_{n+1}+1)(y_{n}+1)}{(x_{n+1}-1)(y_{n}-1)}}=
\frac{1}{1+v_n}\\

v_{n+1}=
\frac{\ln \frac{y_{n+1}+1}{y_{n+1}-1}}{\ln \frac{(y_{n+1}+1)(x_{n}+1)}{(y_{n+1}-1)(x_{n}-1)}}=
\frac{1}{1+u_n} \label{eq8ab}\ea \en

One can solve \eqref{eq8ab}, substitute the solutions in \eqref{eq7b} \& \eqref{eq8b}. In this way we reduce the order of the original system ,,,

Let us now consider the generator $\X_3$. The resulting invariants are 

\bn v_n = \frac{\left[\frac{x_n -1}{x_n +1}\right]^{\alpha_n}}{\frac{y_{n+1}-1}{y_{n+1}+1}}\quad , u_n = \frac{\left[\frac{y_n -1}{y_n +1}\right]^{\beta_n}}{\frac{x_{n+1}-1}{x_{n+1}+1}} \label{eq9b}\en

where \bn \alpha_n = \frac{(1+\sqrt{5})^n-(-1+\sqrt{5})^n}{2[(1+\sqrt{5})^{n-1}+(-1+\sqrt{5})^{n-1}]}, \quad  \beta_n = \frac{(1+\sqrt{5})^n+(-1+\sqrt{5})^n}{2[(1+\sqrt{5})^{n-1}-(-1+\sqrt{5})^{n-1}]} \label{eq9bb}\en

Note also the relationship between them

$$\alpha_{n+1}-1=\frac{1}{\beta_n}, \quad \beta_{n+1}-1=\frac{1}{\alpha_n}$$

From \eqref{eq9b}, we have 

\bn\ba  v_{n+1}= \frac{\left[\frac{x_{n+1}-1}{x_{n+1}+1}\right]^{\alpha_{n}-1}}{\frac{y_n-1}{y_n+1}}=\frac{1}{u_n^{\frac{1}{\beta_n}}},\\

u_{n+1}= \frac{\left[\frac{y_{n+1}-1}{y_{n+1}+1}\right]^{\beta_{n}-1}}{\frac{x_n-1}{x_n+1}}=\frac{1}{v_n^{\frac{1}{\alpha_n}}}

\label{eq10b}\ea \en

One can readily check that the general solution to \eqref{eq10b} is given by

\bn\ba u_n=\frac{1+(-1)^n}{2}\left[ u_0^{\frac{1}{\prod\limits _{k=0}^{\frac{n-2}{2}} \beta_{2k}\alpha_{2k+1}}} \right] + \frac{1-(-1)^n}{2}\left[ v_0^{\frac{-1}{\prod\limits_{k=0}^{\frac{n-1}{2}}\alpha_{2k}\prod\limits_{k=0}^{\frac{n-3}{2}}\beta_{2k+1}}}\right] \\

v_n=\frac{1+(-1)^n}{2}\left[ v_0^{\frac{1}{\prod\limits _{k=0}^{\frac{n-2}{2}} \alpha_{2k}\beta_{2k+1}}} \right] + \frac{1-(-1)^n}{2}\left[ u_0^{\frac{-1}{\prod\limits_{k=0}^{\frac{n-1}{2}}\beta_{2k}\prod\limits_{k=0}^{\frac{n-3}{2}}\alpha_{2k+1}}}\right]
\label{eq11b}\ea\en
where $\alpha_n $ and $\beta_n$ are defined in \eqref{eq9bb}.

The substitution of  \eqref{eq11b} in \eqref{eq9b} reduces the order of the system,,,

\section{Conservation Laws}

In Section 2, we have defined a  first integral  associated to a second-oreder S$\Delta$Es. It is given by \eqref{eq7}

\bn \phi(n,x_n,y_n,x_{n+1},y_{n+1})= \phi(n+1,x_{n+1},y_{n+1},\omega_1,\omega_2)\label{eqa1}\en

Let 

$ P_1=\frac{\p \phi }{\p x_{n}} ,\quad P_2=\frac{\p \phi}{\p x_{n+1}} ,\quad 
\quad Q_1=\frac{\p \phi}{\p y_{n}} ,\quad  Q_2=\frac{\p \phi}{\p y_{n+1}} ,\quad$

By differentiating \eqref{eqa1} with respect to $x_{n},y_{n}, x_{n+1}$ and $y_{n+1}$ we obtain

\bn\ba &P_1=\S(P_2) \omega_{1,x_n} +\S(Q_2)\omega_{2,x_{n}}\\
&Q_1=\S(P_2) \omega_{1,y_n} +\S(Q_2)\omega_{2,y_{n}}\label{eqa2}\ea \en

and \bn\ba&P_2=\S(P_1)+\S(P_2) \omega_{1,x_{n+1}} +\S(Q_2)\omega_{2,x_{n+1}}\\
&Q_2=\S(Q_1)+\S(P_2) \omega_{1,y_{n+1}} +\S(Q_2)\omega_{2,y_{n+1}}\label{eqa3} \ea\en

The substitution of  \eqref{eqa2} in \eqref{eqa3} leads to the following second-order system of functional equations

\bn\ba \S^2(P_2)\S(\omega_{1,x_{n}})+\S^2(Q_2)\S(\omega_{2,x_n})+\S(P_2)
\omega_{1,x_{n+1}}+\S (Q_2)\omega_{2,x_{n+1}}-P_2=0\\
\S^2(P_2)\S(\omega_{1,y_{n}})+\S^2(Q_2)\S(\omega_{2,y_n})+\S(P_2)
\omega_{1,y_{n+1}}+\S (Q_2)\omega_{2,y_{n+1}}-Q_2=0\label{eqa4}
\ea\en

As for SLSC, we differentiate repeatedly to obtain a system of DEs for $P_2$ and $Q_2$.
Given the solutions $P_2,Q_2$ of \eqref{eqa4} we easily construct $P_1,Q_1$. For consistency of our solutions, we must check the integrability conditions

\bn \frac{\p P_1}{\p x_{n+1}}=\frac{\p P_2}{\p x_{n}} \en and \bn 
\frac{\p Q_1}{\p y_{n+1}}=\frac{\p Q_2}{\p y_{n}}\en  
The first integral is then given by

\bn\phi=\int (P_1 \d x_{n}+P_2 \d x_{n+1}+Q_1 \d y_{n}+Q_2 \d y_{n+1}) +F(n)\label{eqa5}\en

The constant of integration $F(n)$ which is a function depending on $n$ is  determined by substituting \eqref{eqa5} in \eqref{eqa1}.

\section{Applications}

Let us consider the second-order S$\Delta$Es

\bn x_{n+2}=a(n) y_n ,\quad y_{n+2}=b(n) x_n \label{eqa6}\en 

By choosing the ansatz $P_2(n,x_n,y_n)$ and $Q_2(n,x_n,y_n)$ one can readily check that   the determining system \eqref{eqa4}  simplifies to

\bn\ba Q_{2}(n+2,\omega_1 ,\omega_2) b(n+1) -P_2(n,x_n,y_n)=0, \\
P_{2}(n+2,\omega_1 ,\omega_2)a(n+1) -Q_2(n,x_n,y_n)=0.
\label{eqa7} \ea \en

where $\omega_1$ and $\omega_2$ denote the right-hand side of \eqref{eqa6}

Differentiating \eqref{eqa7} with respect to $x_n$ and $y_n$ leads to 

\bn P_2=x_n \psi_1(n)+y_n \psi_2(n)+\psi_3(n),\quad Q_{2}= x_n\psi_4 (n)+y_n \psi_5 (n)+\psi_6(n)\label{eqa8}\en 

Thus, we have from \eqref{eqa2}

\bn \ba P_1=b(n)[x_{n+1} \psi_4(n+1) +y_{n+1}\psi_5(n+1)+\psi _6(n+1)]\\

Q_1=a(n)[x_{n+1} \psi_1(n+1) +y_{n+1}\psi_2(n+1)+\psi _3(n+1)] \ea\en

Substituting \eqref{eqa8} in \eqref{eqa7} and separating with respect to $x_n$ and $y_n$ we obtain the system the system 

\bn \ba a(n)b(n+1) \psi_4(n+2)-\psi_2(n)=0,\\
        b(n)a(n+1) \psi_2(n+2) -\psi_4(n)=0,\\
        
       b(n)b(n+1) \psi_5(n+2)-\psi_1(n)=0,\\
        a(n)a(n+1) \psi_1(n+2) -\psi_5(n)=0,\\
        
        b(n+1) \psi_6(n+2)-\psi_3(n)=0,\\
        a(n+1) \psi_3(n+2) -\psi_6(n)=0. \label{eqa9}
  \ea\en
The solutions to  \eqref{eqa9} will provide us the explicit form of $\psi_i, \quad i=1...6$

The first integral is then given by

\bn\ba \phi=& \int b(n)[x_{n+1} \psi_4(n+1) +y_{n+1}\psi_5(n+1)+\psi _6(n+1)]\d x_n +  (x_n \psi_1(n)+y_n \psi_2(n)+\psi_3(n))\d x_{n+1}\\

& a(n)[x_{n+1} \psi_1(n+1) +y_{n+1}\psi_2(n+1)+\psi _3(n+1)] \d y_n 

+(x_n\psi_4 (n)+y_n \psi_5 (n)+\psi_6(n))\d y_{n+1}+K_i

\ea\en
for some constants $K_i$

For clarification, let us consider $a(n)=b(n)=1$, that is,

\bn x_{n+2}= y_n ,\quad y_{n+2}= x_n \label{eqa10}\en 

The solutions to \eqref{eqa9} will be 

\bn\ba \psi_1(n)=\frac{C_1[1+(-1)^n+i^n+(-i)^n]}{4}+\frac{C_2[1-(-1)^n-i(i)^n+i(-i)^n]}{4}+\frac{C_3[1+(-1)^n-i^n-(-i)^n]}{4} +\frac{C_4[1-(-1)^n+i(i)^n-i(-i)^n]}{4}\\
 
 \psi_5(n)=\frac{C_3[1+(-1)^n+i^n+(-i)^n]}{4}+\frac{C_4[1-(-1)^n-i(i)^n+i(-i)^n]}{4}+\frac{C_1[1+(-1)^n-i^n-(-i)^n]}{4} +
 \frac{C_2[1-(-1)^n+i(i)^n-i(-i)^n]}{4}\\

\psi_2(n)=\frac{C_5[1+(-1)^n+i^n+(-i)^n]}{4}+\frac{C_6[1-(-1)^n-i(i)^n+i(-i)^n]}{4}+\frac{C_7[1+(-1)^n-i^n-(-i)^n]}{4} +\frac{C_8[1-(-1)^n+i(i)^n-i(-i)^n]}{4}\\
 
 \psi_4(n)=\frac{C_7[1+(-1)^n+i^n+(-i)^n]}{4}+\frac{C_8[1-(-1)^n-i(i)^n+i(-i)^n]}{4}+\frac{C_5[1+(-1)^n-i^n-(-i)^n]}{4} +
 \frac{C_6[1-(-1)^n+i(i)^n-i(-i)^n]}{4} \\

 \psi_3(n)=\frac{C_9[1+(-1)^n+i^n+(-i)^n]}{4}+\frac{C_{10}[1-(-1)^n-i(i)^n+i(-i)^n]}{4}+\frac{C_{11}[1+(-1)^n-i^n-(-i)^n]}{4} +\frac{C_{12}[1-(-1)^n+i(i)^n-i(-i)^n]}{4}\\
 
 \psi_6(n)=\frac{C_{11}[1+(-1)^n+i^n+(-i)^n]}{4}+\frac{C_{12}[1-(-1)^n-i(i)^n+i(-i)^n]}{4}+\frac{C_9[1+(-1)^n-i^n-(-i)^n]}{4} +
 \frac{C_{10}[1-(-1)^n+i(i)^n-i(-i)^n]}{4}
 \ea\en
where $C_i,\quad i=1...12$ are constants.
We have twelve solutions for $P_2$ and $Q_2$. that is,

1. If $P_2=(\frac{[1+(-1)^n+i^n+(-i)^n]}{4})x_n,\quad   
  Q_2=(\frac{[1+(-1)^n-i^n-(-i)^n]}{4} )y_n$ , then
  
$Q_1=(\frac{[1-(-1)^n+ii^n-i(-i)^n]}{4})x_{n+1},\quad 
 P_1=(\frac{[1-(-1)^n-ii^n+i(-i)^n]}{4} )y_{n+1}$

2. If  $P_2=(\frac{[1-(-1)^n-i(i)^n+i(-i)^n]}{4})x_n, \quad 
  Q_2=(\frac{[1-(-1)^n+i(i)^n-i(-i)^n]}{4})y_n$, then

$Q_1=(\frac{[1+(-1)^n+(i)^n+(-i)^n]}{4})x_{n+1},\quad 
 P_1=( \frac{[1+(-1)^n-(i)^n-(-i)^n]}{4})y_{n+1}$

3. If $P_2=\frac{[1+(-1)^n-i^n-(-i)^n]}{4} )x_n, \quad 
  Q_2=(\frac{[1+(-1)^n+i^n+(-i)^n]}{4})y_n$, then 

$Q_1=(\frac{[1-(-1)^n-ii^n+i(-i)^n]}{4} )x_{n+1},\quad 
 P_1=(\frac{[1-(-1)^n+ii^n-i(-i)^n]}{4})y_{n+1}$

4. If $P_2=(\frac{[1-(-1)^n+i(i)^n-i(-i)^n]}{4})x_n, \quad 
  Q_2=(\frac{[1-(-1)^n-i(i)^n+i(-i)^n]}{4})y_n$, then 

$Q_1=(\frac{[1+(-1)^n-(i)^n-(-i)^n]}{4})x_{n+1},\quad 
 P_1=(\frac{[1+(-1)^n+(i)^n+(-i)^n]}{4})y_{n+1}$

 5. If $P_2=(\frac{[1+(-1)^n+i^n+(-i)^n]}{4})y_n,\quad   
  Q_2=(\frac{[1+(-1)^n-i^n-(-i)^n]}{4} )x_n$ , then
  
$Q_1=(\frac{[1-(-1)^n+ii^n-i(-i)^n]}{4})y_{n+1},\quad 
 P_1=(\frac{[1-(-1)^n-ii^n+i(-i)^n]}{4} )x_{n+1}$

6. If  $P_2=(\frac{[1-(-1)^n-i(i)^n+i(-i)^n]}{4})y_n, \quad 
  Q_2=(\frac{[1-(-1)^n+i(i)^n-i(-i)^n]}{4})x_n$, then

$Q_1=(\frac{[1+(-1)^n+(i)^n+(-i)^n]}{4})y_{n+1},\quad 
 P_1=( \frac{[1+(-1)^n-(i)^n-(-i)^n]}{4})x_{n+1}$

7. If $P_2=\frac{[1+(-1)^n-i^n-(-i)^n]}{4} )y_n, \quad 
  Q_2=(\frac{[1+(-1)^n+i^n+(-i)^n]}{4})x_n$, then 

$Q_1=(\frac{[1-(-1)^n-ii^n+i(-i)^n]}{4} )y_{n+1},\quad 
 P_1=(\frac{[1-(-1)^n+ii^n-i(-i)^n]}{4})x_{n+1}$

8. If $P_2=(\frac{[1-(-1)^n+i(i)^n-i(-i)^n]}{4})y_n, \quad 
  Q_2=(\frac{[1-(-1)^n-i(i)^n+i(-i)^n]}{4})x_n$, then 

$Q_1=(\frac{[1+(-1)^n-(i)^n-(-i)^n]}{4})y_{n+1},\quad 
 P_1=(\frac{[1+(-1)^n+(i)^n+(-i)^n]}{4})x_{n+1}$

9. If $P_2=(\frac{[1+(-1)^n+i^n+(-i)^n]}{4}),\quad   
  Q_2=(\frac{[1+(-1)^n-i^n-(-i)^n]}{4} )$ , then
  
$Q_1=(\frac{[1-(-1)^n+ii^n-i(-i)^n]}{4}),\quad 
 P_1=(\frac{[1-(-1)^n-ii^n+i(-i)^n]}{4} )$

10. If  $P_2=(\frac{[1-(-1)^n-i(i)^n+i(-i)^n]}{4}), \quad 
  Q_2=(\frac{[1-(-1)^n+i(i)^n-i(-i)^n]}{4})$, then

$Q_1=(\frac{[1+(-1)^n+(i)^n+(-i)^n]}{4}),\quad 
 P_1=( \frac{[1+(-1)^n-(i)^n-(-i)^n]}{4})$

11. If $P_2=\frac{[1+(-1)^n-i^n-(-i)^n]}{4} ), \quad 
  Q_2=(\frac{[1+(-1)^n+i^n+(-i)^n]}{4})$, then 

$Q_1=(\frac{[1-(-1)^n-ii^n+i(-i)^n]}{4} ),\quad 
 P_1=(\frac{[1-(-1)^n+ii^n-i(-i)^n]}{4})$

12. If $P_2=(\frac{[1-(-1)^n+i(i)^n-i(-i)^n]}{4}), \quad 
  Q_2=(\frac{[1-(-1)^n-i(i)^n+i(-i)^n]}{4})$, then 

$Q_1=(\frac{[1+(-1)^n-(i)^n-(-i)^n]}{4}),\quad 
 P_1=(\frac{[1+(-1)^n+(i)^n+(-i)^n]}{4})$

%

Therefore, we obtain twelve conservation laws for the system \eqref{eqa10}. They are given by

\bn \ba  \phi_1=& (\frac{[1+(-1)^n+i^n+(-i)^n]}{4})x_n x_{n+1}+
  (\frac{[1+(-1)^n-i^n-(-i)^n]}{4} )y_ny_{n+1}+\\
  &(\frac{[1-(-1)^n+ii^n-i(-i)^n]}{4})x_{n+1}y_n+
 (\frac{[1-(-1)^n-ii^n+i(-i)^n]}{4} )y_{n+1}x_n +K_1\\
 
 \phi_2=& (\frac{[1-(-1)^n-i(i)^n+i(-i)^n]}{4})x_nx_{n+1}+(\frac{[1-(-1)^n+i(i)^n-i(-i)^n]}{4})y_ny_{n+1}+\\
&(\frac{[1+(-1)^n+(i)^n+(-i)^n]}{4})x_{n+1}y_n + 
 ( \frac{[1+(-1)^n-(i)^n-(-i)^n]}{4})y_{n+1}x_n +K_2\\

\phi_3=&\frac{[1+(-1)^n-i^n-(-i)^n]}{4} )x_nx_{n+1} 
  +(\frac{[1+(-1)^n+i^n+(-i)^n]}{4})y_ny_{n+1}+\\

&(\frac{[1-(-1)^n-ii^n+i(-i)^n]}{4} )x_{n+1}y_n 
 +(\frac{[1-(-1)^n+ii^n-i(-i)^n]}{4})y_{n+1}x_n+K_3\\

\phi_4=& (\frac{[1-(-1)^n+i(i)^n-i(-i)^n]}{4})x_nx_{n+1}
 +(\frac{[1-(-1)^n-i(i)^n+i(-i)^n]}{4})y_ny_{n+1}+\\

& (\frac{[1+(-1)^n-(i)^n-(-i)^n]}{4})x_{n+1}y_n 
 +(\frac{[1+(-1)^n+(i)^n+(-i)^n]}{4})y_{n+1}x_n +K_4\\

\phi_5=& (\frac{[1+(-1)^n+i^n+(-i)^n]}{4})y_n x_{n+1}+
  (\frac{[1+(-1)^n-i^n-(-i)^n]}{4} )x_ny_{n+1}+\\
  &(\frac{[1-(-1)^n+ii^n-i(-i)^n]}{4})y_{n+1}y_n+
 (\frac{[1-(-1)^n-ii^n+i(-i)^n]}{4} )x_{n+1}x_n +K_5\\
 
 \phi_6=& (\frac{[1-(-1)^n-i(i)^n+i(-i)^n]}{4})y_nx_{n+1}+(\frac{[1-(-1)^n+i(i)^n-i(-i)^n]}{4})x_ny_{n+1}+\\
&(\frac{[1+(-1)^n+(i)^n+(-i)^n]}{4})y_{n+1}y_n + 
 ( \frac{[1+(-1)^n-(i)^n-(-i)^n]}{4})x_{n+1}x_n +K_6\\

\phi_7=&\frac{[1+(-1)^n-i^n-(-i)^n]}{4} )y_nx_{n+1} 
  +(\frac{[1+(-1)^n+i^n+(-i)^n]}{4})x_ny_{n+1}+\\

&(\frac{[1-(-1)^n-ii^n+i(-i)^n]}{4} )y_{n+1}y_n 
 +(\frac{[1-(-1)^n+ii^n-i(-i)^n]}{4})x_{n+1}x_n+K_7\\

\phi_8=& (\frac{[1-(-1)^n+i(i)^n-i(-i)^n]}{4})y_nx_{n+1}
 +(\frac{[1-(-1)^n-i(i)^n+i(-i)^n]}{4})x_ny_{n+1}+\\

& (\frac{[1+(-1)^n-(i)^n-(-i)^n]}{4})y_{n+1}y_n 
 +(\frac{[1+(-1)^n+(i)^n+(-i)^n]}{4})x_{n+1}x_n +K_8\\

\phi_9=& (\frac{[1+(-1)^n+i^n+(-i)^n]}{4})x_{n+1}+
  (\frac{[1+(-1)^n-i^n-(-i)^n]}{4} )y_{n+1}+\\&
 (\frac{[1-(-1)^n+ii^n-i(-i)^n]}{4})y_n+
 (\frac{[1-(-1)^n-ii^n+i(-i)^n]}{4} )x_n +K_9\\
 
 \phi_{10}=& (\frac{[1-(-1)^n-i(i)^n+i(-i)^n]}{4})x_{n+1}+(\frac{[1-(-1)^n+i(i)^n-i(-i)^n]}{4})y_{n+1}+
\\&(\frac{[1+(-1)^n+(i)^n+(-i)^n]}{4})y_n + 
 ( \frac{[1+(-1)^n-(i)^n-(-i)^n]}{4})x_n +K_{10}\\

\phi_{11}=&\frac{[1+(-1)^n-i^n-(-i)^n]}{4} )x_{n+1} 
  +(\frac{[1+(-1)^n+i^n+(-i)^n]}{4})y_{n+1}+\\

&(\frac{[1-(-1)^n-ii^n+i(-i)^n]}{4} )y_n 
 +(\frac{[1-(-1)^n+ii^n-i(-i)^n]}{4})x_n+K_{11}\\

\phi_{12}=& (\frac{[1-(-1)^n+i(i)^n-i(-i)^n]}{4})x_{n+1}
 +(\frac{[1-(-1)^n-i(i)^n+i(-i)^n]}{4})y_{n+1}+\\

& (\frac{[1+(-1)^n-(i)^n-(-i)^n]}{4})y_n 
 +(\frac{[1+(-1)^n+(i)^n+(-i)^n]}{4})x_n +K_{12}\\
 \ea\en

\end{document}